%% file: sup.tex
\documentclass[12pt]{article}
\usepackage{graphicx}
\usepackage{amsmath,amsfonts,amssymb,amsthm,amsbsy}
\usepackage{fullpage,natbib,verbatim,color}
\usepackage{hyperref,latexsym,mathrsfs}
\input{supmacs}

\newcommand{\eps}{\epsilon}
\renewcommand{\P}{\mathbb{P}}
\renewcommand{\E}{\mathbb{E}}
\newcommand{\Eqn}[1]{(\ref{#1})}
\newcommand{\Thm}[1]{Theorem~\ref{#1}}
\newcommand{\one}[1]{\mathbf{1}_{\{#1\}}}
\newcommand{\Un}{\mbox{Un}}
\newcommand{\Po}{\mbox{Po}}

\definecolor{Purple}        {cmyk}{0.45,0.86,0,0}

\title{On the Supremum of Certain Families of Stochastic Processes}
\author{Wenbo V.  Li \thanks{Department of Mathematical Sciences,
        University of Delaware, \texttt{wli@math.udel.edu}} 
\and    Natesh S.  Pillai\thanks{ CRiSM,
        University of Warwick, \texttt{N.Pillai@warwick.ac.uk}}
\and    Robert L.  Wolpert \thanks{Department of Statistical Science,
        Duke University, \texttt{rlw@stat.duke.edu}}}
\date{\today} 
\begin{document}
\maketitle
\renewcommand{\abstractname}{Summary}
\begin{abstract}
  We consider a family of stochastic processes $\{X_t^\eps, t \in T\}$ on a
  metric space $T$, with a parameter $\eps \downarrow 0$.  We study the
  conditions under which
  \begin{align*}
    \lim_{\eps \to 0} \P \Big( \sup_{ t \in T} |X_t^\eps| < \delta \Big) =1
  \end{align*}
  when one has an \textit{a priori} estimate on the modulus of continuity and
  the value at one point.  We compare our problem to the celebrated
  Kolmogorov continuity criteria for stochastic processes, and finally give
  an application of our main result for stochastic integrals with respect to
  compound Poisson random measures with infinite intensity measures.
\end{abstract}
\section{Introduction}
Let $(T,d)$ be a metric space with finite diameter,
\[
D(T)=\sup \Big \{ d(s,t): s,t \in T \Big \} < \infty.
\]
Let $N(T, d, \delta)$ denote the covering number, \textit{i.e.}, for every
$\delta > 0$, let $N(T, d, \delta)$ denote the minimal number of closed
$d$-balls of radius $\delta$ required to cover $T$.  The supremum of a
stochastic process $X_t$ defined on $T$, $\sup_{t \in T} X_t$ can be
quantified in terms of $N(T,d,\delta)$ (see \cite[Chapter 1]{Tala:05} for
instance) under various assumptions on the process $X_t$.

In this article we consider a family of stochastic processes $X_t^\eps$ on
$T$, with a parameter $\eps > 0$.  In certain applications in nonparametric
statistics (see \Sec{sec:appl}) it is of interest to study the limiting
behaviour of the supremum, $\lim_{\eps \to 0} \sup_{t \in T} X^\eps_t$ when
one has an \textit{a priori} estimate of the form
 \begin{align}
 \E |X_t^\eps -X_s^\eps|^\beta \le B_\eps \,d(s,t)^{\gamma}\notag
 \end{align}
 for some $\beta, \gamma >0$ and $B_\eps \to 0$ as $\eps\to0$.  In
 particular, we would like to identify conditions under which, for every
 $\delta >0$,
 \begin{align}
   \lim_{\eps \to 0} \P \Big( \sup_{ t \in T} |X_t^\eps| < \delta \Big) =1.
   \label{eqn:intmr} 
\end{align}
In our main result in \Sec{sec:mr}, we find conditions in terms of the
covering number $N(T,d,\delta)$ that ensure \eqref{eqn:intmr} holds.
Although our technique is based on well known chaining methods, our principle
result appears to be new.  In \Sec{sec:opt} we discuss briefly the optimality
of our hypotheses and compare our theorem with the Kolmogorov criterion for
continuity of stochastic processes.  In \Sec{sec:appl} we present an
application of our main theorem to random fields constructed from L\'evy
random measures.
%
%

\section{Main Result}\label{sec:mr}
\noindent Let $(T,d)$ be a complete separable metric space and
$(X^\eps_t)_{t\in T}$ a family of real-valued, centered, $L_2$ stochastic
processes on $T$, indexed by $\eps>0$.
   
Let $n_0$ be the largest integer $n$ such that $N(T, d, 2^{-n} ) = 1$ (note
$n_0<0$ is possible).  For every $n \ge n_0$, fix a covering of $T$ of
cardinality $N_n =N(T, d, 2^{-n} )$ by closed balls of radius $2^{-n}$.  From
this we can construct a partition $\cA_n$ of $T$ of cardinality $|\cA_n|=
N_n$ by Borel sets with diameter at most $2^{-n+1}$.  For each $n\ge n_0$,
fix a designated point in each element $A$ of the partition $\cA_n$, and
denote by $T_n$ the collection of these points.  Without loss of generality,
let the designated point in the single element of partition $\cA_{n_0}=\{T\}$
be $T_{n_0}=\{t_0\}$ for a point $t_0\in T$ to be specified in the statement
of \Thm {thm:main} below.
For $t\in T$ denote by $\cA_n(t)$ the partition element $A\in\cA_n$ that
contains $t$.  For every $t$ and every $n$, let $s_n(t)$ be the element of
$T_n$ in $t$'s partition element, so that $t \in \cA_n(s_n(t))$.  It is clear
that $d(t, s_n(t)) \le 2^{-n+1}$ for every $t\in T$ and $n \ge n_0$.  By the
triangle inequality
\begin{align}
  d(s_n(t), s_{n-1}(t))\le d(s_n(t),t) + d(t, s_{n-1}(t)) \le 2^{-n+1}
  +2^{-n+2} = 6 \cdot 2^{-n}. \notag 
\end{align}
Define the set 
\begin{align} 
  H_n \equiv \Big \{(u,v) \in T_n \times T_{n-1}:~ d(u,v)
       \leq 6 \cdot 2^{-n} \Big \}.\label{eqn:setHn} 
\end{align}
The following is our main result: 
\begin{thm}\label{thm:main} Suppose that:
\begin{subequations}
\begin{enumerate}
\item There exists a point $t_0 \in T$ such that
\be
\lim_{\eps \to 0} \E (X_{t_0}^\eps)^2=0.
\label{eqn:var}
\ee
\item There exist $\alpha,\beta>0$ and positive numbers $\{B_\eps\}$ with
  $\lim_{\eps\to0} B_\eps=0$ such that for any $s,t\in T$ \be \E
  |X_t^\eps -X_s^\eps|^\beta \le B_\eps \,d(s,t)^{1+\alpha}.
\label{eqn:Kolm}
\ee
\item There exists a family of partitions $\mathcal{A}_n$ of $T$ of sets of
  diameter no more than $2^{1-n}$ and a constant $\gamma <\alpha$
  such that \be
  \sum_{n=1}^{\infty} |H_n| \,2^{-(1+\gamma)n} < \infty \label{eqn:entrH}
\ee
where $H_n$ is as defined in \eqref{eqn:setHn}.
\end{enumerate}
\end{subequations}
Then for any  $\delta > 0$,
\begin{equation}\label{e:thm1}
  \lim_{\eps \to 0} \P \Big( \sup_{ t \in T} |X_t^\eps| < \delta \Big) =1.
\end{equation}
\end{thm}
\begin{remark}\label{r:kolmo}
  For each fixed $\eps>0$, Equation \eqref{eqn:Kolm} guarantees the existence
  of a path-continuous version of $(X^\eps_t)$ \citep[by Kolmogorov's
  continuity criterion; see][p.\thinspace375.  For more on this connection
  see \Sec{sec:opt}]{Durr:96}.  Since $H_n$ satisfies the bound
\begin{align}
  |H_n| \leq |T_n| \cdot|T_{n-1}| \leq N^2(T, d, 2^{-n}),\notag
\end{align}
the monotonicity of $N(T,d,\delta)$ implies that the entropy condition
\Eqn{eqn:entrH} holds whenever
\begin{align}
  \int_{0}^{D(T)} a^\gamma N^2(T,d,a) \,da < \infty.\notag
\end{align}
Frequently in applications
we have a bound of the form
\begin{align}
  |H_n| \leq C \cdot |T_n| \leq C\cdot N(T,d,2^{-n})\label{eqn:hntn}
\end{align}
for a universal constant $C$ and in this case \eqref{eqn:entrH} holds if
\begin{align}
  \int_{0}^{D(T)} a^\gamma N(T,d,a) \,da < \infty.  \label{eqn:mentrop}
\end{align}
For example, with $T=[0,1]$ and $d(u,v)=|u-v|$, the dyadic partition
\[\cA_n=\left\{\big[i2^{1-n}, (i+1)2^{1-n}\big]:~0\le i<2^{1-n}\right\}\]
of $T$ into $N_n=2^{n-1}$ $d$-balls of radius $2^{-n}$ for $n\ge n_0=1$
satisfies \Eqn{eqn:entrH} for $C=5$.
\end{remark}

\begin{proof}
  Fix $\delta > 0$.  First observe that
\begin{align*}
  \P \Big( \sup_{t \in T} |X_t^\eps| < \delta \Big) 
   &\ge  \P \Big( \sup_{t \in T} |X_t^\e-X_{t_0}^\eps|
             < \delta/2,~ |X_{t_0}^\eps| < \delta/2 \Big) \\
   &\ge \P \Big( \sup_{t \in T} |X_t^\e-X_{t_0}^\eps|
             < \delta/2 \Big)-\P\Big( |X_{t_0}^\eps| \ge \delta/2 \Big)
\intertext{and}
\P( |X_{t_0}^\eps| \ge \delta/2) &\le 4\delta^{-2}\E |X_{t_0}^\eps|^2 \to 0
\end{align*}
as $\eps \to 0$ by equation \eqref{eqn:var}.  Thus we only need to control
$\sup_{t \in T} |X_t^\e-X_{t_0}^\eps|$.

We employ the so-called generic chaining principle of \citet{Ledo:96}
\citetext{see also \citealp{Tala:05}, or \citealp {Xiao:09} for a refinement
  similar in spirit to our approach}.
The fundamental relation is the convergent telescoping sum
\[
X_t - X_{t_0} = \sum_{n>n_0} \left( X_{s_n(t)} -X_{s_{n-1}(t)} \right)
\]
for every $t\in T$, where we note that $s_{n_0}(t)=t_0$ for every $t\in T$.
Then,
\beaa
\sup_{t \in T} | X_t -X_{t_0}|
&\le& \sup_{t \in T} \sum_{n >n_0} | X_{s_n(t)}-X_{s_{n-1}(t)}|  \\
&\le& \sum_{n >n_0} \max_{(u,v) \in H_n} | X_{u}-X_{v}| \eeaa 
For $(u,v) \in T\times T$, let $\{w_n(u,v)\}_{n\geq n_0}$ be a sequence of
non-negative real numbers
such that $\sum_{n \ge n_0} w_n(u,v) = 1$.  For any $\delta >0$, by the
triangle inequality
\begin{align*}
\bigcap_{n > n_0} \bigcap_{(u,v) \in H_n} 
\Big \{| X_u-X_{v}| \leq w_n(u,v) \,\delta/2 \Big \}
\subset \Big \{\sup_{t \in T} | X_t -X_{t_0}| \leq \delta/2 \Big\}.  
\end{align*}
Therefore, 
\beaa \P\left( \sup_{t \in T} |X^\eps_t-X^\eps_{t_0}| > \delta/2
\right)
&\le& \P \Big( \bigcup_{n >n_0} \bigcup_{(u,v) \in T_n}\Big\{
|X^\eps_{u}-X^\eps_{v}| > w_n(u,v) \,\delta/2 \Big\} \Big) \\ 
&\le& \sum_{n >n_0} \sum_{(u,v) \in H_n} \P \left( |X^\eps_{u}-X^\eps_{v}| >
  w_n(u,v)\, \delta/2 \right). \eeaa
Next use Equations \eqref{eqn:Kolm} and \eqref{eqn:entrH} to find optimal
choices for $w_n(v)$ \citep[the so-called ``majorizing measure'', see]
[Chapter 1]{Tala:05}.  Set
\[
  w_n(u,v)  \equiv w_n\equiv (1-2^{-h}) 2^{-h(n-n_0)}, \quad
            h=(\alpha-\gamma)/\beta, ~ v \in T. 
\]
Notice that $\sum_{n \ge n_0} w_n(u,v) = 1$.  By Markov's inequality and
\eqref{eqn:Kolm}, for $v \in T_n$, \beaa \P \left( |X^\eps_{u}-X^\eps_{v}| \ge
  w_n \delta/2 \right)
&\le&  (w_n\delta/2)^{-\beta} \,\E |X^\eps_{u}-X^\eps_{v}|^\beta  \\
&\le&  (\delta/2)^{-\beta}\, (1-2^{-h})^{-\beta}\, 2^{\beta h(n-n_0)}\,
        B_\eps \, d(u,v)^{1+\alpha}  \\ 
&\le& (\delta/2)^{-\beta}\, (1-2^{-h})^{-\beta}\, 2^{-\beta h n_0} \,2^{\beta
  hn}\, B_\eps \, (6 \cdot 2^{-n})^{1+\alpha}. \eeaa
Putting all the estimates together,
\begin{align}
 \P\left( \sup_{t \in T} |X^\eps_t-X^\eps_{t_0}| \ge \delta/2 \right) \notag 
&\le (\delta/2)^{-\beta}\, (1-2^{-h})^{-\beta}\, 6^{1+\alpha}\, 2^{-\beta h
  n_0}\, 
  B_\eps \Big(\sum_{n >n_0} \sum_{(u,v) \in H_n} 2^{\beta hn}\,
  2^{-(1+\alpha)n }\Big) \notag \\ 
&= C  B_\eps \sum_{n >n_0} |H_n|  2^{-(1+\gamma)n } \notag \notag  
\label{eqn:finest}
\end{align}
for a finite constant $C<\infty$.  Since $B_\eps \to 0$ as $\eps \to 0$ and
the sum converges by \eqref{eqn:entrH},
\[
\lim_{\eps \to 0}\P\left( \sup_{t \in T} |X^\eps_t-X^\eps_{t_0}| \ge
  \delta/2 \right) = 0
 \] 
and the theorem is proved.
\end{proof}

\section{Near Optimality of Our Hypothesis}\label{sec:opt}

Kolmogorov's continuity criterion asserts the existence of a path continuous
version of any stochastic process $X_t, t \in [0,1]$ that satisfies
\begin{align}
 \E\Big(|X_t - X_s|^\beta \Big) \leq C~|t-s|^{1+\alpha} \notag 
\end{align}
for some fixed $\alpha, \beta >0$, $C<\infty$ and all $0\le s,t\le 1$ (cf.
\Eqn{eqn:Kolm}).  Strict inequality $\alpha>0$ is necessary, as illustrated
by the well known counter example
\begin{align*}
  X_t &= \one{U \le t}
\end{align*}
for $U \sim \Un[0,1]$ which satisfies $\E |X_t - X_s|^\beta \le C |t-s|$ for
all $\beta>0$ and $C\ge1$ but is almost surely discontinuous.  In the spirit
of this example, here we construct a stochastic process which shows that our
hypothesis (2) in \Thm{thm:main} is ``very close'' to optimal.

Let $U \sim \Un[0,1]$, $0<\eps<1$ and $X_t^\eps=\one{t < U \le t+\eps}$, $0
\le t \le 1$.  Then for any fixed $t \in [0,1]$,
\begin{align}
  \E (X_t^\eps)^2&=\P(t< U \le t+\eps)=\min(\eps, 1-t)\label{e:eg.l1}
\intertext{Since}
\E X_t^\eps X_s^\eps &=
  \begin{cases}
    0 & \text{if } |t-s| > \eps \\
    \e-|t-s| &\text{if }  |t-s| \le \eps,~ \min(s,t) \le 1-\eps \\
    1-\max(s,t) &\text{if } |t-s| \le \eps,~ \min(s,t) \ge 1-\eps \\
  \end{cases},\notag\\
\intertext{it follows that}
  \E (X_t^\e-X_s^\eps)^2 &\le 2 \min(\eps , |t-s|).\label{eqn:examest}
\end{align}
By \Eqn{e:eg.l1} $X_t^\eps$ satisfies \Eqn{eqn:var} for any $t_0\in[0,1]$,
and by \Eqn{eqn:examest} we have bounds on $\E[(X_t^\eps - X_s^\eps)^\beta]$
for $\beta=2$ both of the form $B|t-s|$ (for fixed $B=2$) and of the form
$B_\eps\to0$ (with $B_\eps=2\eps$), but not quite a bound of the form
required by \Eqn{eqn:Kolm}.  The conclusion \Eqn{e:thm1} of \Thm{thm:main}
fails for the process $X^\eps_t$ since, for any $\eps >0$, $\sup_{0 \le t \le
  1}X_t^\eps=1$ almost surely.
We believe that the condition $\alpha >0$ in equation \eqref{eqn:Kolm} cannot
be relaxed and state this a conjecture.
\begin{conj}
  \Thm{thm:main} is not true if hypothesis (2) (equation
  \eqref{eqn:Kolm}) is replaced by
\begin{align*}
  \E |X_t^\eps -X_s^\eps|^\beta \le B_\eps \,d(s,t).
\end{align*}
\end{conj}
%
%
\section{Application: Compensated Poisson Random Measures}\label{sec:appl}
In this section we present an application of \Thm{thm:main} to a stochastic
process constructed from compensated Poisson random measures.

Let $\O$ be a Polish space and $\nu(du\,d\o)$ be a positive
sigma-finite measure on $(-1,1)\times \O$ such that
\begin{align*}
\nu\big((-a,a)\times \O\big) &= \infty, ~~ \forall a \in [0,1] \\
\iint_{(-1,1)\times\Omega}  u^2~\nu(du\,d\o)&< \infty.
\end{align*}
Let
\begin{equation}
N(du\,d\o) \sim \Po(\nu)\notag
\end{equation}
be a Poisson random measure on $(-1,1)\times\O$ which assigns independent
$\Po\big(\nu(B_i)\big)$ distributions to disjoint Borel sets $B_i \subset
(-1,1)\times\O$.  Let
\begin{equation}
\tilde{N}(du\,d\o) \equiv N(du\,d\o) - \nu(du\,d\o)\notag
\end{equation}
denote the compensated Poisson measure with mean $0$, an isometry from
$L_2\big((-1,1)\times\O,\nu(du\,d\o)\big)$ to the square-integrable
zero-mean random variables \citep[\textit{p.} 38]{Sato:99}.

Let $K(t,\o):[0,1]\times \O \to \bbR$ be a Borel measurable
function such that
\begin{align}
  \iint_{(-1,1)\times \O} K^2(t,\o)~u^2~\nu(du\,d\o) < \infty \label{eqn:intcond}
\end{align}
for all $0\le t\le1$.  For $0<\eps\le1$ define a stochastic process
$X^\eps_t$ by
\begin{align}
  X^\eps_t \equiv \iint_{\{0 < |u| <\eps\}\times \O}
  K(t,\o)\,u\,\tilde{N}(du\,d\o), ~~ t \in [0,1].\label{eqn:condstoc}
\end{align}
For every $t \in [0,1]$ the stochastic integral \eqref{eqn:condstoc} is well
defined by \eqref{eqn:intcond} \citep[see][]{Wolp:Taqq:05,Rajp:Rosi:89}.  For
$t \in [0,1]$, we have:
\begin{subequations}
\begin{align}
  \E\left[X_t^\eps \right] &= 0 \notag\\
  \E\left[(X_t^\eps) ^2 \right]&= \iint_{(-\eps,\eps)\times \O}
  K^2(t,\o)~u^2~\nu(du\,d\o) < \infty \notag\\ 
  \E\left[ e^{i\zeta X_t^\eps} \right] &=
    \exp \Big \{ \iint_{\mathbb(-\eps, \eps)\times \O}
    \left[ e^{i\zeta K(t,\o)u}-1 - i\zeta K(t,\o)u\right]\nu(du\,d\o)
    \Big\},\notag 
\end{align}
\end{subequations}
the L\'evy-Khinchine formula for the characteristic function of an
infinitely divisible random variable.

The stochastic process $\big \{X^\eps \equiv X_t^\eps, t \in [0,1] \big \}$
is the discretization error arising from the approximation of certain
stochastic integrals by finite sums \citep[see][]{Pill:Wolp:08,
  Wolp:Clyd:Tu:2006}.  The limiting behaviour of the process $X^\eps$ as
$\eps$ goes to zero \citep[see][\S3]{Pill:Wolp:08} is of particular interest;
we would like to identify the conditions on the function $K$ under
which\begin{equation} \lim_{\eps \to 0} \P\Big(\sup_{t \in [0,1]}|X_t^\eps| >
  \delta \Big) = 0
\label{eqn:concpois}
\end{equation}
for all $\delta > 0$, so the approximation error vanishes in the limit.
Concentration equalities similar to \eqref{eqn:concpois} were studied by
\citet{Bour:06} for finite intensity measures (\textit{i.e.,}
$\nu((-1,1)\times \O) < \infty$) using methods that are not applicable to our
infinite intensity case.

In the next proposition we apply \Thm{thm:main} to identify conditions for
the kernel $K(\cdot,\cdot)$ under which \eqref{eqn:concpois} holds.
\begin{prop}\label{p:example}
  Let $K(t,\o):[0,1]\times \O \to \bbR$ satisfy \Eqn{eqn:intcond} and
\begin{equation}
  |K(t,\o) - K(s,\o)|^2 \leq C(\omega) ~|t - s|^{1+ \alpha}, \quad s,t \in
  [0,1],~\omega\in\Omega \label{eqn:klip}
\end{equation}
for some $\alpha >0$ and Borel measurable function $C: \Omega \to \bbR_+$ satisfying
\begin{equation}
  \iint_{(-1,1)\times \O} C(\omega) ~u^2 ~\nu(du\,d\o) < \infty.\label{eqn:cwint}
\end{equation}
Let $X^\eps$ be the stochastic process on $[0,1]$ given in \eqref{eqn:condstoc}.
Then, for any $\delta > 0$,
\begin{align}
  \lim_{\eps \to 0} \P\Big(\sup_{t \in [0,1]}|X_t^\eps| > \delta \Big) = 0.\notag
\end{align}
\end{prop}
\begin{proof}
 For any $t_0 \in [0,1]$, 
by \Eqn{eqn:intcond} and the dominated convergence theorem
\begin{align}
  \lim_{\eps \to 0}\E\Big[\big( X^\eps_{t_0} \big)^2\Big] = \lim_{\eps \to
    0}\int_{\{|u| \le \eps\}\times \O} K^2(t_0,\o)~u^2~\nu(du\,d\o) = 0,
\label{eqn:verih1}
\end{align}
verifying hypothesis 1 (equation \eqref{eqn:var}) of \Thm{thm:main}.  For
$t,s \in [0,1]$, by \eqref{eqn:klip} and by the isometric property of
$\tilde{N}(du\,d\o)$,
\begin{align}
  \E\Big[\big( X^\eps_t - X^\eps_s\big)^2\Big]
  &=  \iint_{\{|u| \leq \eps\}\times \O}
      |K(t,\o) - K(s,\o)|^2~u^2~\nu(du\,d\o)\notag\\ 
  &\leq B_\eps |t - s|^{1+ \alpha},\text{\qquad where}\notag\\
  B_\eps &\equiv \iint_{\{|u| \leq \eps\}\times \O}
  C(\omega)~u^2~\nu(du\,d\o)\to0 \label{eqn:diffest}
\end{align}
as $\eps\to0$ by \eqref{eqn:cwint}, so hypothesis 2 (equation \ref{eqn:Kolm})
of \Thm{thm:main} is satisfied with the Euclidean metric $d(t,s) \equiv |t -
s|$.  For dyadic partitions of $[0,1]$, we have already shown that
\Eqn{eqn:hntn} holds.  Since $N([0,1],d,a) =\lceil\frac1a\rceil \leq
\frac{2}{a}$ for all $0<a<1$, for any $\gamma>0$ (say, $\gamma =\alpha/2$),
\begin{align}
\int_{0}^1 a^\gamma N([0,1],d,a)\,da \leq 2\int_{0}^1 a^{\gamma - 1} da 
  =\frac2\gamma< \infty
\label{eqn:covver}
\end{align}
verifying Equation \eqref{eqn:mentrop}.  Therefore by
\eqref{eqn:verih1},\eqref{eqn:diffest},\eqref{eqn:covver} and \Thm{thm:main},
it follows that for any $\delta >0$,
\[
\lim_{\eps \to 0} \P\Big(\sup_{t \in
  [0,1]}|X_t^\eps| > \delta \Big) = 0
\]
and we are done.
\end{proof}

\begin{remark}
  It is not known whether the conclusion of the above proposition still holds
  if \eqref{eqn:klip} is weakened to
\begin{align}
  |K(t,\o) - K(s,\o)|^2 &\leq C(\omega) ~|t - s|, ~s,t \in [0,1].\notag
\end{align}
\end{remark}
\section*{Acknowledgments}
The authors thank the anonymous referee and Prof.\ Yimin Xiao for their
thoughtful comments.  This work was supported in part by the National Science
Foundation under Grant Numbers
DMS-0805929, 
DMS-0757549, 
DMS-0635449  
and the CRiSM research fellowship.  Any opinions, findings, and conclusions
or recommendations expressed in this material are those of the authors and do
not necessarily reflect the views of the National Science Foundation.
\bibliographystyle{authordate1} 
\bibliography{statjour-abbr,sup}
\end{document}

%% file: supmacs.tex
\newtheorem{thm}{Theorem}

\newtheorem{prop}{Proposition}

\newtheorem{conj}{Conjecture}

\newcommand{\be}{\begin{equation}}
\newcommand{\ee}{\end{equation}}
\newcommand{\bea}{\begin{eqnarray}}
\newcommand{\eea}{\end{eqnarray}}
\newcommand{\beaa}{\begin{eqnarray*}}
\newcommand{\eeaa}{\end{eqnarray*}}
\newcommand{\bei}{\begin{itemize}}
\newcommand{\eei}{\end{itemize}}
\newcommand{\bee}{\begin{enumerate}}
\newcommand{\eee}{\end{enumerate}}
\newcommand{\bi}{\begin{itemize}}
\newcommand{\ei}{\end{itemize}}
\newcommand{\beq}{\begin{eqnarray*}}
\newcommand{\eeq}{\end{eqnarray*}}
\newcommand{\beqn}{\begin{eqnarray}}
\newcommand{\eeqn}{\end{eqnarray}}

\newcommand{\ignore}[1]{}{}

\newcommand{\e}{\varepsilon}

\def\<{\langle}
\def\>{\rangle}

\def\P{{\mathbb P}}

\def\E{{\mathbb E}\,}

\def\O{\Omega}
\def \o{\omega}

\newcommand{\Sec}[1]{Section~\ref{#1}}

\newcommand{\cA}{{\cal A}}

\newcommand{\bbR}{\mathbb{R}}

%% file: sup.bbl
\begin{thebibliography}{}

\bibitem[\protect\citename{Durrett, }1996]{Durr:96}
Durrett, Richard. 1996.
\newblock {\em Probability: theory and examples}. Second edn.
\newblock Belmont, CA: Duxbury Press.

\bibitem[\protect\citename{Ledoux, }1996]{Ledo:96}
Ledoux, Michel. 1996.
\newblock Isoperimetry and {G}aussian analysis.
\newblock {\em Pages  165--294 of:} {\em Lectures on probability theory and
  statistics (Saint-Flour, 1994)}.
\newblock Lecture Notes in Mathematics, vol. 1648.
\newblock Berlin, DE: Springer-Verlag.

\bibitem[\protect\citename{Pillai \& Wolpert, }2008]{Pill:Wolp:08}
Pillai, Natesh~S., \& Wolpert, Robert~L. 2008.
\newblock {\em Posterior Consistency of {B}ayesian Nonparametric Models Using
  {L\'e}vy Random Field Priors}.
\newblock Discussion Paper 2008-08. Duke Univ. Dept. Statist. Science, USA.

\bibitem[\protect\citename{Rajput \& Rosi{\'n}ski, }1989]{Rajp:Rosi:89}
Rajput, Balram~S., \& Rosi{\'n}ski, Jan. 1989.
\newblock Spectral representations of infinitely divisible processes.
\newblock {\em Probab. Theory Rel.}, {\bf 82}(3), 451--487.

\bibitem[\protect\citename{Reynaud-Bouret, }2006]{Bour:06}
Reynaud-Bouret, Patricia. 2006.
\newblock Compensator and exponential inequalities for some suprema of counting
  processes.
\newblock {\em Statist. Probab. Lett.}, {\bf 76}(14), 1514--1521.

\bibitem[\protect\citename{Sato, }1999]{Sato:99}
Sato, Ken-iti. 1999.
\newblock {\em L{\'e}vy Processes and Infinitely Divisible Distributions}.
\newblock Cambridge Studies in Advanced Mathematics, vol. 68.
\newblock Cambridge, UK: Cambridge Univ. Press.

\bibitem[\protect\citename{Talagrand, }2005]{Tala:05}
Talagrand, Michel. 2005.
\newblock {\em The generic {C}haining}.
\newblock Springer Monographs in Mathematics.
\newblock Berlin, DE: Springer-Verlag.

\bibitem[\protect\citename{Wolpert \& Taqqu, }2005]{Wolp:Taqq:05}
Wolpert, Robert~L., \& Taqqu, Murad~S. 2005.
\newblock Fractional {O}rnstein-{U}hlenbeck {L\'e}vy Processes and the
  {T}elecom Process: Upstairs and Downstairs.
\newblock {\em Signal Processing}, {\bf 85}(8), 1523--1545.

\bibitem[\protect\citename{Wolpert {\em et~al.\ }\relax,
  }2006]{Wolp:Clyd:Tu:2006}
Wolpert, Robert~L., Clyde, Merlise~A., \& Tu, Chong. 2006.
\newblock {\em L{\'e}vy {A}daptive {R}egression {K}ernels}.
\newblock Discussion Paper 2006-08. Duke Univ. Dept. Statist. Science, USA.
\newblock Revised November 2009.

\bibitem[\protect\citename{Xiao, }2009]{Xiao:09}
Xiao, Yimin. 2009.
\newblock On Uniform Modulus of Continuity of Random Fields.
\newblock {\em Monatshefte f{\"u}r Mathematik}.
\newblock To appear.

\end{thebibliography}
